\documentclass[12pt, reqno]{amsart}
\usepackage[utf8]{inputenc}
\usepackage{amsmath,amsthm, color, amscd, amsfonts}
\usepackage{amsfonts}
\usepackage{amssymb}
\usepackage{graphicx}
\usepackage{hyperref}
\newtheorem{thm}{Theorem}[section]
 \newtheorem{cor}[thm]{Corollary}
 \newtheorem{lemma}[thm]{Lemma}
 \newtheorem{prop}[thm]{Proposition}
 \theoremstyle{definition}
 
 \theoremstyle{remark}
 \newtheorem{rem}[thm]{Remark}

 \numberwithin{equation}{section}
 \usepackage{hyperref}
 \graphicspath{{/texfiles/}}
  \usepackage{amssymb}

\begin{document}
\title[Unitary subgroup and normal complement problem]{On the Normal Complement Problem for Finite Group Algebras of Abelian-by-Cyclic Groups}


\author{Allen Herman }
\address{Department of Mathematics and Statistics \\ University of Regina \\ Saskatchewan \\ Canada} \email{allen.herman@uregina.ca}
\author{Surinder Kaur }
\address{Department of Mathematics \\ School of Engineering and Sciences \\ SRM University-AP \\ Andhra Pradesh, India-522502} \email{skaur.math@gmail.com}



\begin{abstract}
Assume $F$ is a finite field of order $p^f$ and $q$ is an odd prime for which $p^f-1=sq^m$, where $m \ge 1$ and $(s,q)=1$. In this article, we obtain the order of symmetric and unitary subgroup of the semisimple group algebra $FC_q.$ Further, for the extension $G$ of $C_q = \langle b \rangle$ by an abelian group $A$ of order $p^n$ with $C_{A}(b) = \{e\}$, we prove that if $m>1,$ or $(s+1) \geq q$ and $2n \geq f(q-1)$, then $G$ does not have a normal complement in $V(FG)$.
\end{abstract}

\maketitle

\section{Introduction}\label{sec1}

In this paper we will consider the normal complement problem for modular group algebras $FG$, when $F$ is a finite field of characteristic $p$ and $G \simeq A \rtimes C_q$ is a non-abelian semi-direct product of a finite abelian $p$-group $A$ and a cyclic group of prime order $q \ne p$. The normal complement problem originates in Graham Higman's thesis (see \cite{MR0002137}), where it is shown that $G$ has a torsion-free normal complement in $\mathcal{U}(\mathbb{Z}G)$ when $G$ is a finite abelian group. In \cite{ward1960some}, H. N. Ward posed the normal complement problem for $FG$ in the case when $G$ is a finite $p$-group and $F$ is the field with $p$ elements; that is, whether or not $\mathcal{U}(FG)$ has a normal subgroup $N$ for which $\mathcal{U}(FG) = N \rtimes G$. In \cite{MR0444697}, R. K. Dennis asked for which groups $G$ and commutative rings $R$ does $G$ have a normal complement in $\mathcal{U}(RG)$? Since that time the normal complement problem has received quite a lot of attention in research on group rings - to bring readers up-to-date on these results we refer them to \cite{MR457539}, \cite{MR447403}, \cite{MR560574}, \cite{MR1045598}, \cite{MR3590872},  \cite{MR3969485},  \cite{MR4115375},  \cite{MR4442491},  \cite{MR4542679}, \cite{MR1969234}, \cite{MR1231619}, \cite{MR1221725}, \cite[Chapter 4]{MR1242557} and \cite{MR4370430}. For the groups we will deal with in this paper, existence of a torsion-free normal complement to $G$ in $\mathcal{U}(\mathbb{Z}G)$ follows from a result of Cliff, Sehgal, and Weiss \cite{MR0641639} (see also \cite[(31.1b)]{MR1242557}). This is in contrast to our main result, which shows our groups $G$ will not have a normal complement in $\mathcal{U}(FG)$ when $|A|$ and the $q'$-part of $|F^\times|$ are large enough.

This paper is divided into four sections. In Section 2, we introduce our notation and give some preliminary results on the unit group of $FG$. In Section 3, we study the conjugacy classes of some units of $FC_q$ in the group of normalized units $V(FG).$ In Section $4,$ we prove the main theorem of the paper:

\begin{thm} \label{mainthm}
Let $p$ be an odd prime and $F$ be the field with $p^f$ elements. Assume that $q$ is an odd prime divisor of $(p^f-1)$ such that $p^f = (sq^m+1),$ where $(s,q)=1$.
Let $G$ be the extension of the cyclic group $B = \langle b \rangle$ of order $q$ by an abelian group $A$ of order $p^n$ such that $C_A(b) = \{e\}$. If $m>1,$ or $(s+1) \geq q$ and $2n \geq f(q-1)$, then $G$ does not have a normal complement in $V(FG).$
\end{thm}

The case of $m>1$ was established earlier in \cite{MR3590872} (see Theorem 2). Thus this is a new result in the case $m=1.$ 


\section{Notation and Preliminaries}
Throughout the paper, $p$ will be a prime, $F$ will be a finite field of characteristic $p$, $G$ will be a finite group, $FG$ will be the corresponding group algebra, and $\mathcal{U}(FG)$ will be its group of units. The augmentation ideal $Aug(FG)$ is the kernel of the augmentation map $\epsilon: FG \rightarrow F$ defined by the rule $\displaystyle \sum_{g \in G} \alpha_g g \mapsto \sum_{g \in G}\alpha_g.$ The units of augmentation 1 form a subgroup which is called the normalized unit
group of $FG$ and is denoted by $V(FG).$ Note that $\displaystyle \mathcal{U}(FG) =
F^{\times} \times V(FG)$ and hence the existence of a
normal complement of $G$ in $\displaystyle \mathcal{U}(FG)$
implies the existence of a normal complement of $G$ in 
$\displaystyle V(FG)$ and vice-versa.

The inverse map $G \to G$ extends F-linearly to an anti-automorphism $*: FG \to FG$ as $\displaystyle \Big(\sum_{g \in G} \alpha_g g \Big)^{*} = \sum_{g \in G} \alpha_g g^{-1}.$ Elements $k \in FG$ are said to be symmetric if $\displaystyle k^{*} = k$ and skew-symmetric if $k^{*} = -k.$ Furthermore, the units in $V(FG)$ inverted by $\displaystyle *$ form a group which is called the unitary subgroup and is denoted by $\displaystyle V_{*}(FG).$ A study of the symmetric and unitary units of the modular group algebras of finite $p$-groups has been done in \cite{MR1283327}, \cite{MR1872239}, \cite{MR1747361}, \cite{MR1374635} and \cite{MR3974958}. In \cite{MR1853809}, the authors studied the unitary subgroup of the group algebra of a finite group over a non-absolute field of odd characteristic. 
Note that in particular, if $G$ is an abelian group, then the symmetric units form a group and we denote the symmetric subgroup of $V(FG)$ by $V_{+}(FG)$. Further, observe that if $x \in V_{*}(FG) \cap V_{+}(FG),$ then $x^{*} = x = x^{-1}$ implies that $x^2=1.$ Thus the subgroup $ V_{*}(FG) \cap V_{+}(FG)$ is an elementary abelian $2$-group.

For a normal subgroup $H$ of $G,$ we denote the two-sided ideal of $FG$ generated by the
elements $\displaystyle \{ h-1 \ | \ h \in H \}$ by $\displaystyle \Gamma(H).$ 
For a transversal $T$ of $H$ in $G,$ the elements $\displaystyle \{ t(h-1) 
\ | \ t \in T, h \in H \}$ form a basis of $\displaystyle \Gamma(H)$ as an $F$-module. If $H$ is a finite $p$-group, then $Aug(FH)$ is a nil ideal of $FH$, and $V(FH)$ is equal to $1 + Aug(FH)$. If $H$ is a normal $p$-subgroup of $G,$ then $\displaystyle \Gamma(H) = Aug(FH).FG$ is a nilpotent ideal of $FG$ and so $\displaystyle (1+\Gamma(H))$ is a normal subgroup of $\displaystyle V(FG).$

Let $\displaystyle G = A \rtimes B,$ where $A$ is a finite $p$-group. Then the map $\displaystyle \rho:G \rightarrow B$ defined by the rule $a \mapsto e$ and $b \mapsto b,$ where $a \in A$ and $b \in B,$ can be extended linearly to $F$-algebra homomorphism from $FG$ to $FB$ whose kernel is $\displaystyle \Gamma(A).$ If $\rho'$ denotes the group homomorphism induced by $\rho$ on $V(FG),$ then $\rho'$ induces the isomorphism \[\frac{V(FG)}{(1+\Gamma(A))} \cong V(FB).\]
If $\displaystyle i: V(FB) \rightarrow V(FG)$ denotes the natural inclusion, then $\displaystyle \rho' \circ i$ is an identity map on $V(FB)$ and hence the short exact sequence \[\{1\} \rightarrow (1+\Gamma(A)) \rightarrow V(FG)  \xrightarrow{\rho'} V(FB)  \rightarrow \{1\}\] splits. Therefore, 
$\displaystyle V(FG) = (1+\Gamma(A)) \rtimes V(FB).$ Note that if $B$ is an abelian group, then $\displaystyle V'(FG) \subseteq (1+\Gamma(A))$ where $\displaystyle V'(FG)$ denotes the derived subgroup of $V(FG).$

Let $K$ be a finite group and $p$ be a prime divisor of $|K|$. Assume that $p^t$ is the maximal power of $p$ dividing $|K|.$ If $K$ is solvable, then it has a Hall $p'$-subgroup, a subgroup of index $p^t.$ We denote this subgroup by $K_{p^{'}}$. For any element $x$ in a group $K,$ we denote the centralizer of $x$ in $K$ by $C_K(x)$ and the conjugacy class of $x$ in $K$ by $Cl_x.$ For any integer $n,$ let $C_n$ denote the cyclic group of order $n.$

Our first result gives information about the normalized unit group $V(FB)$ when $|F|=p^f > 3$ and $B$ is cyclic of order $q$ for some odd prime $q$ dividing $p^f-1.$
In the following result, which is a generalization of Proposition $3.1$ from \cite{MR4542679}, we obtain the order of the subgroups $V_{*}(FB)$ and $V_{+}(FB)$ and establish a structure decomposition of $V(FB)$. Assume that $\zeta$ is a primitive element of $F^{\times}$ and $\omega =\zeta^{(p^f-1)/q}$ is a root of unity of order $q$. Then we have:

\begin{prop}\label{prop1}
Let $F$ be a field with $p^f > 3$ elements and $q$ be an odd prime divisor of $(p^f-1)$.
Assume that $B = \langle b \rangle$ is the cyclic group of order $q.$ Then \[ |V_{+}(FB)| = |V_{*}(FB)| = |F^{\times}|^{\frac{q-1}{2}}\] and \[ (V(FB))_{2^{'}} = (V_{+}(FB))_{2^{'}} \times (V_{*}(FB))_{2^{'}}.\]

\end{prop}

\begin{proof}
It follows from Chinese remainder theorem that \[FB \cong \frac{F[x]}{\langle x^q-1 \rangle} \cong \bigoplus_{i=0}^{(q-1)} \frac{F[x]}{\langle (x-\zeta^i) \rangle}\cong  \underbrace{F \oplus F \oplus \ldots \oplus F}_ {q \ copies}.\] Here note that the map $\displaystyle \psi: FB \rightarrow F \oplus F \oplus \ldots \oplus F$ defined by the rule $\displaystyle b \mapsto (1, \zeta, \zeta^2, \ldots, \zeta^{q-1})$ is an isomorphism. Therefore, 
\begin{equation} \label{eqn1}
 V(FB) \cong \underbrace{F^{\times} \times F^{\times} \times \ldots \times F^{\times}}_{(q-1)  \ copies}.
 \end{equation}
 
Since $FB$ is commutative and semisimple, one can obtain its basis constituted by the primitive idempotents by using the character formula: 

$e_0 = \frac{1}{q}(1+b+ \dots + b^{q-1})$, 

$e_1 = \frac{1}{q}(1 + \omega^{q-1} b + \omega^{q-2} b^2 + \dots + \omega b^{q-1}),$

$e_2 = \frac{1}{q}(1 + \omega^{q-2}b + \omega^{q-4}b^2 + \dots + \omega^2 b^{q-1}),$

$\qquad \vdots$ 

$e_{q-1} = \frac{1}{q}(1 + \omega b + \omega^2b^2 + \dots + \omega^{q-1}b^{q-1})$. 

\smallskip

Therefore, any element $u \in FB$ can be expressed as $\displaystyle \sum_{i=0}^{q-1} ue_i = \sum_{i=0}^{q-1} u_i e_i,$ where $u_i \in F$. In particular, if $u$ is a unit of $FB,$ then it has the property that $ue_i \ne 0$ for all $0 \leq i \leq (q-1)$. 

Note that the $F$-linear involution $b^* = b^{-1}$ of $FB$ fixes $e_0$ and reverses the list $e_1,e_2,\dots,e_{q-2},e_{q-1}$. Thus $\displaystyle u^* = u_0e_0 + \sum_{i>0} u_{q-i} e_i$. If $u$ is a unitary unit, then it satisfies $u^*=u^{-1}$ and so $u_{q-i} = u_i^{-1}$. 
Similarly, if $u$ is a symmetric unit, then it follows from $u^*=u$ that $u_{q-i}=u_i$. Therefore, $\displaystyle |V_{*}(FB)|=|V_{+}(FB)|=(p^f-1)^{(q-1)/2}$ and hence it follows from the isomorphism in \ref{eqn1} that each of the subgroups $V_{*}(FB)$ and $V_{+}(FB)$ is isomorphic to $C_{p^f-1}^{(q-1)/2}$. 

It follows from the discussion in earlier part of this section that the intersection of $V_+(FB)$ and $V_*(FB)$ is an elementary abelian $2$-group. Therefore, the Hall $2^{'}$-subgroups of  $V_+(FB)$ and $V_*(FB)$ intersect trivially.   
\end{proof}

In the next Lemma, we will use the previous Proposition to identify unitary units $u \in V(FB)$ for which $F[u] = FB$:

\begin{lemma} \label{span}
Suppose $u \in \mathcal{U}(FB)$ has $q$ distinct projections $ue_i = u_ie_i,$ where
$u_i \in F.$ Then the $F$-span of $\langle u \rangle$ is all of $FB.$ In particular, there exists a
polynomial $p(x)$ of degree at most $(q-1)$ for which $b = p(u).$
\begin{proof}

Since $FB$ is a commutative semisimple algebra, the set of primitive idempotents $e_0, \dots, e_{q-1}$ introduced in Proposition \ref{prop1} forms a basis of $FB$. Considering the regular representation $FB \rightarrow M_q(F)$ written in terms of the primitive idempotent basis, if
$$ u = u_0 e_0 + u_1 e_1 + \dots + u_{q-1}e_{q-1} $$
for $u_i \in F$ satisfying $u e_i = u_i e_i$ for $i=0,1,\dots,q-1$, then the matrix $[u]$ in this representation will be the diagonal matrix with diagonal $(u_0,u_1,\dots,u_{q-1})$.

Writing elements of $FB$ as column vectors in terms of the same basis produces the regular module for the regular representation.  Since $u e_i = u_i e_i$, the $e_i$ correspond to eigenvectors of $[u]$   with eigenvalues $u_i$ in this module, for every $i=0,\dots,q-1$. It follows that $[u]$ has $q$ linearly independent eigenvectors and is hence diagonalizable.  Since $[u]$ is a diagonalizable matrix, the degree of its minimal polynomial is equal to the number of distinct eigenvalues $u_i$ of $[u]$. The degree of the minimal polynomial of $[u]$ is the dimension of the subring generated by the matrix $[u]$ in $M_q(F)$.  Since the regular representation is faithful, this will also be the dimension of the subalgebra of $FB$ generated by $u$. So from this discussion we can conclude that any $u \in FB$ having $q$ distinct eigenvalues $u_i$ will have the property that $F[u] = FB$. The Lemma now follows.
\end{proof}
    
\end{lemma}

\begin{rem}

If $u \in V(FB)$ is a unitary unit of order $q$ with $q$ distinct eigenvalues, then the map $b \mapsto u$ induces an augmentation-preserving $F$-linear $*$-automorphism of $FB$.  For example, in $F_{11} C_5$, the unitary unit
$$ u = e_0 + \omega e_1 + \omega^3 e_2 + \omega^2 e_3 + \omega^4 e_4 = 2b + 3b^2 +8b^3 +10 b^4 $$
is such a unit.  For this unit $u$, the polynomial $p(x) = 2 x + 3 x^2 + 8 x^3 + 10 x^4$ has the curious feature that both $u=p(b)$ and $b = p(u)$! So in this case $b \mapsto u$ is an automorphism of order $2$ that swaps polynomials in $b$ with polynomials in $u$, and vice-versa. 
\end{rem}

\section{Conjugacy Classes of some Unitary Units of $FC_q$ in $V(F(A \rtimes C_q))$}

Let $\displaystyle G = A \rtimes B$ be a non-abelian group, where $A$ is a finite $p$-group and $B$ is a cyclic group of odd order $q$. It follows from the discussion in Section $2$ that $\displaystyle V(FG) = (1+\Gamma(A)) \rtimes V(FB).$ 
Since $V(FB)$ is abelian, for any element $\displaystyle w$ of $V(FB), C_{V(FG)}(w) = C_{(1+\Gamma(A))}(w) \rtimes V(FB)$. Thus we have the following:

\begin{rem}\label{rem1}
Let $w$ be an element of $V(FB).$ Then the length of the conjugacy class of $w$ in $V(FG)$ is $\displaystyle |Cl_w| = \frac{|V(FG)|}{|C_{V(FG)}(w)|} = \frac{|(1+\Gamma(A))|}{|C_{(1+\Gamma(A))}(w)|}.$
\end{rem}

In the next section, we solve the normal complement problem for $V(FG)$ and the elements of $V(FB)$ with the largest class length in $V(FG)$ will play a crucial role towards this end. Note that it follows from the above remark that $\displaystyle |Cl_b| = \frac{|(1+\Gamma(A))|}{|C_{(1+\Gamma(A))}(b)|}.$ Since $C_{(1+\Gamma(A))}(b) \subseteq C_{(1+\Gamma(A))}(w)$, the elements of $V(FB)$ whose class length in $V(FG)$ is the largest are the ones whose centralizer in $V(FG)$ (or in $(1+\Gamma(A)))$ is same as that of $b.$  

Now the following corollary to Lemma \ref{span} provides a sufficient condition for an element of $V(FB)$ to have its centralizer in $V(FG)$ be the same as that of $b$:

\begin{cor} \label{centb}
If $u \in \mathcal{U}(FB)$ has $q$ distinct projections $ue_i = u_ie_i,$ where $u_i \in F,$ then $C_{V(FG)}(u) = C_{V(FG)}(b).$ Thus for such element $u, C_{(1+\Gamma(A))}(u) = C_{(1+\Gamma(A))}(b).$
\end{cor}
\begin{proof}
In Lemma \ref{span}, we identify unitary units $u \in V(FB)$ for which $F[u] = FB$. Of course, when $F[u] = FB,$ then $b \in F[u]$, which implies $C_{V(FG)}(u) \subseteq C_{V(FG)}(b)$, and hence these are equal.\end{proof}

Note that it follows from Proposition \ref{prop1} that the Sylow $q$-subgroup of $V_{*}(FB)$ is isomorphic to the direct product of $\displaystyle \frac{(q-1)}{2}$ copies of $C_q.$ Then since $B$ is a subgroup of $V_{*}(FB),$ one can write $V_{*}(FB) = W \times B,$ where $W$ is a normal complement of $B$ in $V_{*}(FB).$ When  $q>3$, $W$ is isomorphic to the direct product of $\Big( \frac{(q-1)}{2}-1 \Big)$ copies of $C_q.$ 

In particular, assume that the characteristic $p$ of $F$ is such that $|F^{\times}| = (p^f-1) = sq^m,$ where $q>3$ and $(s,q)=1.$ Then as an application of Corollary \ref{centb}, we obtain in the following that if $(s+1) \ge q$, then any normal complement of $B$ in the unitary subgroup $V_{*}(FB)$ contains an element whose centralizer in $V(FG)$ is same as that of $b$.

\begin{lemma} \label{centsameb} 
Let $(p^f-1) = sq^m,$ where $q>3$ and $(s,q)=1.$ Then
\begin{enumerate}
\item If $m>1,$ B has no complement in $V_*(FB).$ 
\item If $m = 1$ and $(s+1) \ge q,$ then every complement of $B$ in $V_*(FB)$ contains a unitary unit $u$ whose centralizer in $V(FG)$ is the same as that of $b.$ 
\end{enumerate}
\end{lemma}
\begin{proof}

(1): Since every subgroup of order $q$ in $V_*(FB)$ is contained in a cyclic subgroup of order $q^m,$ it follows that $B$ cannot have a complement in $V_*(FB).$ We (remark) here that the case of $m>1$ was also established in \cite[Theorem 2]{MR3590872}.

(2): 
Owing to Corollary \ref{centb}, it is sufficient to show that every complement of $B$ in $V(FB)$ contains a unitary unit with $q$ distinct projections. Let $N$ be a complement of $B$ in $V(FB).$ Clearly, $N = N_{q^{'}} \times N_q.$ Suppose there is an element $n \in N_q$ that does not have $q$ distinct projections. But since $n$ is a unitary and is of order $q$, it has at least two projections that are not equal to $1.$ Without loss of generality, one can assume these are $n_1$ and $n_{q-1}=n_1^{-1}.$ 

Now we aim to obtain a unit $v$ in $N_{q'}$ such that $v_0 = v_1 = v_{q-1} = 1$; and that the rest of the projections $\{ v_i \ | \ i=2,\dots,(q-2)\}$ are distinct. The desired conclusion would then follow from Corollary \ref{centb} because now the unitary unit $w=vn$ is an element of $N$ with distinct projections. Let $\eta$ be a primitive $s$-th root of unity in $F^\times$. The desired aim can be fulfilled if we choose the $(q-3)$ elements $\{ v_i \ | \ i=2,\dots,(q-2)\}$ from distinct powers of $\eta$ in such a way that $v_{q-i} = v_{i}^{-1}$ while making sure that the choice is not $1$ or $-1$. In other words, the $(q-3)$ positions have to be filled with the remaining $(s-2)$ powers of $\eta$ which can be done successfully if and only if $(s-2) \ge q-3.$  This proves (ii).
 
\end{proof}

In the upcoming discussion, we assume that $A$ is abelian of order $p^n$. We first prove that in this case, $C_{(1+\Gamma(A))}(b)$ is an abelian group. Indeed, assume that $\displaystyle \alpha = 1+\sum_{a \in A} \alpha_a (a-1) \in C_{(1+\Gamma(A))}(b),$ where $\alpha_a \in FB$. If the orbits of $A$ under the action of $B$ are $O_{a_1}, O_{a_2}, \ldots, O_{a_l}$ with the representatives $a_1, a_2, \ldots a_l$, then $b^{-1} \alpha b = \alpha$ implies that $\displaystyle \alpha = 1+\sum_{i=1}^{l} \alpha_{a_i} (\widehat{O_{a_i}}-|O_{a_i}|),$ where $\widehat{O_{a_i}}$ denotes the sum of the elements of orbit $O_{a_i}.$ Clearly, such elements commute with each other and hence the claim follows.

 
Now we list the following results on the length of conjugacy classes of $V(FG)$ which go a long way in resolving normal complement problem for $FG$. We assume that exponent of $(1+\Gamma(A))$ is $p^k.$ The following lemma is the first result in the aforementioned direction and is a generalization of \cite[Lemma 2]{MR4115375}:

\begin{lemma}\label{lem3}
Let $S= \{w \in V(FB) \ | \ C_{(1+\Gamma(A))}(w) = C_{(1+\Gamma(A))}(b)\}$. Then the elements of the set $U = \{wz \ | \ w \in S, z \in C_{(1+\Gamma(A))}(b)\}$ lie in disjoint conjugacy classes of $V(FG).$
\end{lemma}
\begin{proof}
Assume that $w_1z_1, w_2z_2 \in U$ are such that $v^{-1}(w_1z_1)v = w_2z_2,$ for some $v \in V(FG)$. As $V'(FG) \subseteq (1+\Gamma(A)),$ it follows that $z_1(w_1z_1,v)z_{2}^{-1} = w_1^{-1}w_2 \in (1+\Gamma(A)) \cap V(FB) = \{1\}$. Thus $w_1 = w_2.$ Now observe that $\displaystyle v^{-1}(w_1z_1)^{p^k}v = (w_2z_2)^{p^k}$ and $v^{-1}(w_1z_1)^{|F^{\times}|}v = (w_2z_2)^{|F^{\times}|}$ imply $v^{-1}w_1v = w_2 = w_1$ and $v^{-1}z_1v = z_2$. Therefore, $v \in C_{(1+\Gamma(A))}(w_1)$. Further, since $w_1 \in S, C_{V(FG)}(w_1) = C_{(1+\Gamma(A))}(b) \rtimes V(FB)$ and so we can write $v = \alpha \beta,$ for some $\alpha \in C_{(1+\Gamma(A))}(b)$ and $\beta \in V(FB).$ As $C_{(1+\Gamma(A))}(b)$ is an abelian group, $v^{-1}z_1v = z_2$ implies $\beta^{-1}z_1\beta = z_2$ which gives $z_1=z_2.$
\end{proof}

The following lemma, which is a consequence of Lemma 3 in \cite{MR4115375}, provides a lower bound on the class length of elements of $U$:

\begin{lemma} \label{lowerbound}
Assume that $\gamma$ is an element of $V(FB)$. If $\displaystyle z \in C_{(1+\Gamma(A))}(\gamma)$, then $|Cl_{\gamma z}| \geq |Cl_{\gamma}|$.
\end{lemma}

\begin{proof}
Let $v \in C_{(1+\Gamma(A))}(\gamma z).$ Then it follows from $v^{-1}(\gamma z)^{p^k}v=(\gamma z)^{p^k}$ and $v^{-1}(\gamma z)^{|F^\times|}v=(\gamma z)^{|F^\times|}$ that $v^{-1}\gamma v=\gamma$ and $v^{-1}zv=z$. Therefore, we have $C_{(1+\Gamma(A))}(\gamma z) \subseteq C_{(1+\Gamma(A))}(\gamma) \cap C_{(1+\Gamma(A))}(z).$ Since the reverse inclusion is obvious, we get that 
$C_{(1+\Gamma(A))}(\gamma z)  = C_{(1+\Gamma(A))}(\gamma) \cap C_{(1+\Gamma(A))}(z).$ Now the result follows.
\end{proof}

\section{The Unitary Subgroup and the Normal Complement Problem in $V(F(A \rtimes C_q))$}
In this section, we employ the information about the conjugacy classes in the unit group of $FG$ to answer the normal complement problem. The following is a proof of the main result of this paper:

\vspace{.2cm}

\textbf{Proof of Theorem 1.1:}
First note that when $q=3,$ the result follows from Theorem 2 of \cite{MR4115375}. Therefore, we now prove it for any prime $q>3.$ Suppose $G$ has a normal complement, say $N,$ in $V(FG).$ Then there exists an epimorphism $\displaystyle \delta: V(FG) \rightarrow G$ which fixes the elements of $G$ and that $\displaystyle \ker(\delta) = N.$ Thus $\displaystyle V(FG) = N \rtimes G.$ The restriction of the map $\delta$ to $\displaystyle V(FB)$ is also an epimorphism and we denote its kernel $\displaystyle N \cap V(FB)$ by $N(FB)$. Indeed, if for any $\displaystyle w \in V(FB), \delta(w) = ab^j,$ where $a \in A$ and $0 \leq j \leq (q-1),$ then $\displaystyle \delta(w b^{q-j}) = a,$ which implies that $a = e.$ Therefore, the existence of a normal complement of $G$ in $V(FG)$ implies the existence of a normal complement of $B$ in $V(FB)$. Note that it follows from Lemma \ref{centsameb} that if $m>1$, then $B$ does not have a normal complement in $V(FB)$ and hence $G$ can not have one in $V(FG).$ Therefore, we assume that $m=1$. 

Now we give the following lemma which reduces the normal complement problem of $G$ in $V(FG)$ to that of $G$ in $\displaystyle V_{*}(FG):$
\begin{lemma} \label{lemuni}
If $G$ has a normal complement in $V(FG),$ then $G$ has one in $\displaystyle V_{*}(FG).$
\end{lemma}
\begin{proof}
The restriction of $\delta$ to $\displaystyle V_{*}(FG)$ implies that $\displaystyle N_{*}(FG) = N \cap V_{*}(FG)$ is a normal complement of $G$ in $\displaystyle V_{*}(FG).$
\end{proof}

Now note that $\displaystyle V_{*}(FG) = (1+\Gamma(A))_{*} \rtimes V_{*}(FB),$ where $(1+\Gamma(A))_{*}$ is subgroup of unitary units of $(1+\Gamma(A))$. Thus $\displaystyle V_{*}(FG) = (1+\Gamma(A))_{*} \rtimes (N_{*}(FB) \times B),$ where $ \displaystyle N_{*}(FB) = N \cap V_{*}(FB)$ is a normal complement of $B$ in $\displaystyle V_{*}(FB).$

Further, assume that $S_1$ and $S_2$ denote the set of symmetric and skew-symmetric elements in $\displaystyle \Gamma(A).$ Since $p$ is an odd prime, any element $\displaystyle v \in \Gamma(A)$ can be expressed as $\displaystyle \Big( \frac{v+v^{*}}{2} \Big) + \Big( \frac{v-v^{*}}{2} \Big) \in S_1 \oplus  S_2.$ Therefore, $\displaystyle \Gamma(A) = S_1 \oplus S_2.$ 

Consequently, we obtain the following relation between the size of the conjugacy class of a unitary unit of $V(FB)$ in $V(FG)$ and that in $V_{*}(FG)$: 

\begin{lemma} \label{lemsqrt}
Let $x$ be a unitary unit in $FB.$ Then \[|C_{(1+\Gamma(A))_{*}}(x)| = |C_{(1+\Gamma(A))}(x)|^{\frac{1}{2}}.\]
\end{lemma}

\begin{proof}
Note that if $\displaystyle x \in V(FB)$ is unitary, then for any $\displaystyle y \in C_{\Gamma(A)}(x),$ $\displaystyle y^{*} \in C_{\Gamma(A)}(x).$ Indeed $\displaystyle xy = yx$ implies $\displaystyle y^{*}x^{*} = x^{*}y^{*}$ and hence $\displaystyle y^{*}x = xy^{*}.$ As $\displaystyle y = \frac{y+y^{*}}{2} + \frac{y-y^{*}}{2},$ $\displaystyle C_{\Gamma(A)}(x) = C_{S_1}(x) \oplus C_{S_2}(x).$ Clearly, there is a bijection between the elements of $C_{S_1}(x)$ and $C_{S_2}(x).$ Thus $\displaystyle |C_{\Gamma(A)}(x)| = |C_{S_1}(x)|^2 = |C_{S_2}(x)|^2$. Note that it follows from Lemma $9$ in \cite{MR4115375} that each skew-symmetric element $l$ in $\displaystyle \Gamma(A)$ corresponds to the Cayley unitary unit $\displaystyle (1-l)(1+l)^{-1}$ in $\displaystyle (1+\Gamma(A))_{*}$ and vice-versa. Hence $\displaystyle |(1+\Gamma(A))_{*}| = |S_2|;$ and so $\displaystyle |C_{(1+\Gamma(A))_{*}}(x)| = |C_{S_2}(x)|.$ Therefore, $\displaystyle |C_{(1+\Gamma(A))_{*}}(x)|^2 = |C_{\Gamma(A)}(x)|.$ 
\end{proof}

If $\displaystyle Cl^{*}_{x}$ denotes the conjugacy class of unitary unit $x$ of $V(FB)$ in $\displaystyle V_{*}(FG),$ then the above discussion implies that \[ |Cl_{x}| = \frac{|(1+\Gamma(A))|}{|C_{(1+\Gamma(A))}(x)|} = \frac{|(1+\Gamma(A))_{*}|^2}{|C_{(1+\Gamma(A))_{*}}(x)|^2} = |Cl^{*}_{x}|^2.\]

We remark that Lemma \ref{lem3} holds for the conjugacy classes restricted to the unitary subgroup as well. We have:

\begin{lemma}\label{starlem}
Let $T= \{w \in V_{*}(FB) \ | \ C_{(1+\Gamma(A))_{*}}(w) = C_{(1+\Gamma(A))_{*}}(b)\}$ and $W = \{wz \ | \ w \in T, z \in C_{(1+\Gamma(A))_{*}}(b)\}.$ Then the elements of $W,$ when conjugated by elements of $V_{*}(FG),$ lie in disjoint conjugacy classes.
\end{lemma}

Similarly, we have the following result in which we rephrase Lemma \ref{lowerbound} restricted to $V_{*}(FG)$:

\begin{lemma} \label{starlowerbound}
Let $\gamma$ be an element of $V_{*}(FB)$. If $\displaystyle z \in C_{(1+\Gamma(A))_{*}}(\gamma)$, then $|Cl_{\gamma z}^{*}| \geq |Cl_{\gamma}^{*}|$.
\end{lemma}



If $B$ has a normal complement in $V(FB),$ then it follows from the discussion succeeding Lemma \ref{lemuni} that $\displaystyle V_{*}(FB) = N_{*}(FB) \times B.$ Observe that for any $\displaystyle z \in C_{(1+\Gamma(A))}(b),$ $\displaystyle z^{-1}bz = b$ implies that $\displaystyle \delta(z^{-1})b \delta(z) = b.$ Therefore, $\displaystyle \delta(z) \in C_A(b).$ Since $\displaystyle C_A(b)$ is trivial, $\displaystyle C_{(1+\Gamma(A))}(b) \subseteq N.$

Recall that it follows from the proof of Lemma \ref{centsameb} that if the prime $p$ is such that $|F^{\times}| = p^f-1 = sq^m$ with $(s+1) \geq q,$ then $N_{*}(FB)$ contains an element $w$ whose centralizer in $V(FG)$ is same as that of $b$. Further, it is obvious from the way the element $w$ is chosen that its order is a multiple of $q$. Therefore, if we let $R$ denote the set $\{w^i \ | \ (i, o(w)) = 1 \},$ then $|R| > \phi(q) = (q-1).$

Now consider the set $\displaystyle T^{'} = \{ xz \ | \ x \in R, z \in C_{(1+\Gamma(A))_{*}}(b)\}.$ According to Lemma \ref{starlem}, the elements of $\displaystyle T^{'}$ are representatives of disjoint conjugacy classes in $V_{*}(FG).$ Here observe that all the elements of $\displaystyle T^{'}$ and their conjugacy classes lie in $\displaystyle N_{*}(FG) \setminus (1+\Gamma(A))_{*}.$ Further, it follows from Lemma \ref{starlowerbound} that $\displaystyle |Cl_{xz}^{*}| \geq |Cl_x^{*}|,$ where $\displaystyle |Cl_{x}^{*}| = |Cl_b^{*}| = \frac{|(1+\Gamma(A))_{*}|}{|C_{(1+\Gamma(A))_{*}}(b)|}.$
Therefore, the number of elements in the conjugacy classes of elements of $T^{'}$, \[ | \bigcup_{t \in T^{'}} Cl_t^{*} | > (q-1) |C_{(1+\Gamma(A))_{*}}(b)| |Cl_{x}^{*}| = (q-1) |(1+\Gamma(A))_{*}|.\]
 
On the other hand, 
\begin{align*}
|N_{*}(FG) \setminus (1+\Gamma(A))_{*}| & = |N_{*}(FG)| - |N_{*}(FG) \cap (1+\Gamma(A))_{*}| 
\\ & = \frac{|V_{*}(FG)|}{|G|} - \frac{|(1+\Gamma(A))_{*}|}{|A|}
\\ & = \frac{|(1+\Gamma(A))_{*}|}{|A|} \Big( \frac{|V_{*}(FB)|}{|B|} - 1 \Big).
\\ & = \frac{|(1+\Gamma(A))_{*}|}{|A|} \Big( \frac{(sq)^{\frac{q-1}{2}}}{q} - 1 \Big).
\end{align*}
As $p = (sq+1)^{\frac{1}{f}},$ $|A| = (sq+1)^{\frac{n}{f}} > \Big( \frac{(sq)^{\frac{n}{f}}}{q}-1\Big).$ Further, since $\frac{n}{f} \geq \frac{(q-1)}{2},$ it follows that $|A|$ exceeds $\Big( \frac{(sq)^{\frac{q-1}{2}}}{q} - 1 \Big).$ Thus 

\begin{align*}
|\bigcup_{t \in T^{'}} Cl_t^{*}| & > (q-1) | 
(1+\Gamma(A))_{*}| > | 
(1+\Gamma(A))_{*}| = |A| \frac{| 
(1+\Gamma(A))_{*}|}{|A|} \\ & > \Big( \frac{(sq)^{\frac{q-1}{2}}}
{q} - 1 \Big)\frac{|(1+\Gamma(A))_{*}|}{|A|} \\ 
& = |N_{*}(FG) \setminus (1+\Gamma(A))_{*}|,
\end{align*}

which is a contradiction.

 
\qed

\begin{cor}
Let $p$ be a prime, where $p = (s5^m+1) > 11$ with $(s,5)=1$ and $F$ be the field with $p$ elements. Consider the extension $\displaystyle G$ of $C_5 = \langle b \rangle$ by an abelian group $A$ of order $p^n$ with $n \geq 2$ such that $C_{A}(b) = \{e\}$. Then $G$ does not have a normal complement in $V(FG).$
\end{cor}
\begin{proof}
When $m>1,$ it follows from Lemma \ref{centsameb} that $B$ has no complement in $V_{*}(FB)$ and hence $G$ can not have one in $V_{*}(FG).$ Thus we only need to consider the case when $m=1.$ Since $p > 11$, we have that $s \geq 6$ and so $(s+1) > 5.$ Further, if $A$ is an abelian group of order $p^n$ with $n \geq 2 = \frac{(q-1)}{2}$, then it follows from Theorem \ref{mainthm} that $G$ does not have a normal complement in $V(FG).$ 

 \end{proof}

\section{Acknowledgement}
The authors are grateful to the anonymous referee for the constructive feedback which improved the overall presentation of the article and led to the modification of the main result. SK thanks the research support of the National Board for Higher Mathematics, Department of Atomic Energy, Govt. of India (Grant Number: 02011/1/2025/NBHM(R.P.)/R\&D II/1128)

\end{document}